\input amstex
\documentstyle{amsppt} \magnification=1000
\hcorrection{0mm} \vcorrection{0mm} \hsize=15.8cm \vsize=23cm
\pageno=1 \NoBlackBoxes \nologo \nopagenumbers 

\def\blr#1#2{ { \buildrel{#2} \over {#1} } }
\def\trm#1{ {\text{\rm{#1}}} }
\def\ra{\rightarrow} \def\lra{\longrightarrow} \def\ms{\mapsto}
\def\dim{{\operatorname{dim}}}
\def\image{{\operatorname{image}}}
\def\ker{{\operatorname{ker}}} 
\def\Aut{{\operatorname{Aut}}\,}
\def\smooth{{\operatorname{smooth}}} \def\sing{{\operatorname{sing}}}
\def\hom{{\operatorname{hom}}} \def\alg{{\operatorname{alg}}}

\def\OO{{\Cal O}} \def\ZZ{{\Cal Z}} \def\DD{{\Cal D}} \def\JJ{{\Cal J}} \def\HH{{\Cal H}}
\def\Z{\Bbb Z} \def\Q{\Bbb Q}  \def\C{\Bbb C} \def\P{\Bbb P}

\def\prg#1#2{\vskip5mm\parskip=1mm\parindent=4mm\noindent{{\bf{\S \ #1.\qquad #2.}}}\vskip3mm}

\def\Thm#1{\noindent {\bf Theorem #1.}\ \it}
\def\rmk#1{\noindent{\bf Remark #1. \ }}

\def\proof#1{\vskip1mm\parindent=0mm{\it\underbar{Proof #1}. \ }\parskip=0mm}
\def\qed{ \vskip-1pt {\rightline{$\square$}} \parskip=1mm \parindent=4mm }

 at 8truept

\def\zfst{1} \def\zfnf{2} \def\zfnr{3} \def\gtha{4} \def\rmkd{5} \def\gthb{6} \def\gthbp{}



\topmatter

\title
A Griffiths' Theorem for varieties with isolated singularities
\endtitle

\author
Vincenzo Di Gennaro, \ Davide Franco, \ Giambattista Marini
\endauthor

\leftheadtext{Vincenzo Di Gennaro, \ Davide Franco, \ Giambattista Marini}
\rightheadtext{A Griffiths' Theorem for varieties with isolated singularities}

\abstract
By the fundamental work of Griffiths one knows that, under suitable assumption, homological and algebraic equivalence do not
coincide for a general hypersurface section of a smooth projective variety $Y$. In the present paper we
prove the same result in case $Y$ has isolated singularities.


\endabstract
\endtopmatter

\document
\baselineskip 12pt
\parskip=1mm \parindent=4mm

\prg{0}{Introduction}

By the fundamental work of Griffiths [Gri] it is well known that
homological equivalence and algebraic equivalence do not coincide
in higher codimension. Griffiths' results rest on a careful
analysis of the normal function associated to an algebraic and
primitive cycle in a smooth projective hypersurface $ Y $. One of
the main steps of his analysis consists in showing that the
algebraic part of the intermediate Jacobian of a general
hyperplane section of $ Y $ vanishes, because of the
irreducibility of the monodromy action on the rational cohomology.

A deep improvement of Griffiths' results was made by Nori (see
[N]) who proved that homological and algebraic equivalence do not
coincide also in a range where the Jacobian fibrations are trivial
by homological reasons, so that it seems hard to deduce Nori's
results from properties of normal functions. Moreover Nori's work
shows that Griffiths' result holds not only for algebraically
trivial cycles but for all the cycles induced by correspondences
from homologically trivial cycles on other varieties. All the
analysis follows from the celebrated Nori's Connectedness Theorem,
inspired by previous works of Green [Gre] and Voisin
(unpublished). A complete  account of this beautiful story can be
found e.g. in [V2].

In this paper we investigate the same questions for varieties with isolated singularities
(see Theorems \gtha \ and \gthb \ below for precise statements).
Unfortunately, as far as we know, it is not known a connectedness theorem for singular
varieties so our analysis is more classical and strongly based on the use of normal functions.
One of our main ingredients is a monodromy theorem for varieties with isolated singularities
proved recently in [DGF]. Most of the steps of the proofs are non
trivial for singular varieties so we found convenient to add a lot
of details which usually do not appear in the literature.

\prg{1}{Notation and preliminaries}

Throughout this paper $ Y \subseteq \P^N $ denotes an irreducible, reduced, projective variety
having at worst isolated singularities, with
$$
\dim\,Y \quad = \quad n+1 \quad = \quad 2r \quad \ge \quad 4 \ .
$$
Furthermore
\ $ f : \, \widetilde Y \ra Y \, $ denotes a resolution of singularities;
\ $ \Sigma_d \subseteq \vert \OO_Y(d) \vert $ denotes the linear system on $ Y $ cut by $ d $-degree hypersurface sections,
\ $ \Sigma_d^\circ := \left\{ b \in \Sigma_d \big\vert X_b \trm{ is smooth} \right\} $
\ its subset parameterizing smooth varieties (for $ b \in \Sigma_d $ $ X_b $ denotes the corresponding hypersurface section)
and \ $ \DD \, := \, \Sigma_d \setminus \Sigma_d^\circ $ \ the discriminant locus;
\ for $ \ b \, \in \, \Sigma_d^\circ \ $,
$ \ i_b \, : \ X_b \, \ra \, \widetilde Y $ \ denotes the natural inclusion;
\ for a analytical submanifold $ \ B \ $ of $ \ \Sigma_d^\circ \ $ we consider the natural family over $ \, B \, $
$$
Y_B \quad := \quad \big\{ (x,\,b) \ \in \ Y \times B \ \big\vert \ x \in X_b \ \big\} , \qquad \pi \, : \ Y_B \, \ra \, B \, ,
$$
and the natural inclusion $ \ i : \, Y_B \hookrightarrow \widetilde Y \times B \ $ of families over $ \, B \, $ obtained by globalizing $ \, i_b $.
We consider the monodromy representation associated to the universal family, i.e. over $ \Sigma_d^\circ $,
$ \ \pi_{{}_1}(\Sigma_d^\circ,\, b) \, \ra \, \Aut \big( H^n(X; \, \Q) \big) $
where $ \ b \, \in \, \Sigma_d^\circ \ $ is any point and $ \, X = X_b \, $ ([V2], Ch. 3).
We want to recall that as a consequence of Deligne invariant subspace Theorem the set $ \, I \, $ of invariant elements is
the image of the pull-back in cohomology $ \ i_b^* \, : \ H^n(\widetilde Y; \, \Q) \, \ra \, H^n(X; \, \Q) \, $,
in particular it is a Hodge substructure of $ \ H^n(X; \, \Q) \ $.
Its orthogonal complement shall be denoted by $ \, V \, $ and will be called the \lq\lq vanishing cohomology". So, $ \ H^n(X;\Q) \, = \, I \perp V \, $.

We also consider the analytic map of Jacobian fibrations over $ \, B \, $ associated
to the inclusion $ \, Y_B \subseteq \widetilde Y \times B \, $,
which we shall denote by $ \, i^* \, $
$$
\CD
J^r(\widetilde Y) \, \times \, B @. \qquad @> i^* >> \qquad @. \JJ @. \qquad := \qquad \big\{ J^r(X_b) \big\}_{b\in B}.
\endCD
$$
This map globalizes the map of Griffiths' intermediate Jacobians $ \ i_b^* \, : \ J^r(\widetilde Y) \, \ra \, J^r(X_b) \ $,
where we keep the standard notation $ \, J^r(W) \, := \, H^{2r-1}(W;\C) \big/ \big( F^r H^{2r-1}(W;\C) \oplus H^{2r-1}(W;\Z) \big) \, $
for the Griffiths' intermediate Jacobian of any smooth projective variety $ \, W \, $.
\ The kernel of such $ \, i_b^* \, $ does not depend on $ \ b \, $, in other terms the inverse image
via $ \, i^* \, $ above of the trivial section $ \, (i^*)^{-1}(\{0\}_{b\in B}) \subseteq \JJ $ \ is a product
$ \, K \times B \, $, where $ \, K \, $ is a subgroup of $ \, J^r(\widetilde Y) \, $.
We furthermore denote with $ \, T \, $ the image subtori-fibration of $ \, J^r(\widetilde Y) \times B \, $ in $ \, \JJ \, $ and with
$ \ T_b \, = \, i_b^* \big( J^r(\widetilde Y) \big) \ $ its fiber over $ \, b \, $.
We have an exact sequence of families over $ \, B \, $
$$
\CD
0 \quad @>>> \quad K \ \times \ B \quad @>>> \quad J^r(\widetilde Y) \, \times \, B \quad @>>> \quad T \quad @>>> \quad 0 \ .
\endCD
$$
The inclusion $ \ T \, \subseteq \, \JJ \ $ is closed and does not depend on the resolution of singularities
$ \widetilde Y \, \ra \, Y \, $, likewise any inclusion $ \, T_b \subseteq J^r(X_b) \, $.

For any $ m $-dimensional projective variety $ W , \ \ZZ^q(W) \, , \ \ZZ^q(W)_{\hom} \, $ and $ \, \ZZ^q(W)_{\alg} \, $
denote respectively the group of $ q $-codimensional algebraic cycles,
its subgroup of homologically trivial cycles and its subgroup of algebraically trivial cycles
(compare with [F1], Ch. 19 and [F2], Appendix B).
For a cycle $ \ Z \in \ZZ^q(W) \ $
we shall denote the corrsponding classes in the $ \, q $-Chow group and in homology respectively as
$ \ [Z] \in CH^{q}(W) \ $ and $ \ cl(Z) \in H_{2m-2q}(W;\,\C) \, $.
The notation for the Abel-Jacobi map (in case $ W $ is smooth) shall be $ \, \Psi_{_{AJ}} : \ZZ^q(W)_{\hom} \to J^q(W) $,
we also set $ J^q(W)_\alg := \Psi_{_{AJ}}\big(\ZZ^q(W)_{\alg}\big) $, the image of algebraically trivial cycles.
Furthermore, in case $ \, W \, $ is a smooth quasi-projetive variety, $ \, cl(Z)^\vee \, $ denotes the
cohomology class of the cycle $ Z $ in $ \, H^{2q}(W;\,\C) \, $ ([V1], Ch. 11).
We also recall that any local complete intersection (l.c.i.) morphism between projective varieties induces Gysin maps in homology and in cohomology
([F1], 19.2.1).

Let us go back to our situation. For a smooth $ \, d $-degree hypersurface section $ \, X \, = \, X_b \, $
we define $ \, J^r(X)_0 \, $ as the cokernel of the map of Jacobians $ \ J^r(\widetilde Y) \, \ra \, J^r(X) $
\ and eventually we set $ \ J^r(X)_{0,\,\alg} $ \ as the image of $ \, J^r(X)_\alg \, $ in $ \, J^r(X)_0 \, $.

\rmk{\zfst}
Let $ \, Z \in \ZZ^{r}(Y) \, $.
If $ \ [Z] \vert_{_X} \, \in \, CH^r(X)_\hom \ $ for some smooth $ \, d $-degree hypersurface section $ \ X \, \hookrightarrow \, Y \, $, \ then
$ \ [Z] \vert_{_{X_b}} \, \in \, CH^r(X_b)_\hom \ $ \ for all smooth $ \, d $-degree hypersurface sections $ \ X_b \, \hookrightarrow \, Y \, $.

To prove this observe that the image of $ cl(Z) $ in $ H^{2r}(X_b;\C) $ via Gysin map
and Poincar\'e Duality
$ H_{2n+2-2r}(Y;\C)\to H_{2n-2r}(X_b;\C)\cong H^{2r}(X_b;\C)$ is equal to the
image of $cl(\widetilde Z)^\vee$ via restriction map $H^{2r}(\widetilde Y;\C)\to H^{2r}(X_b;\C)$,
where $ \widetilde Z \in \ZZ^r (\widetilde Y) $ denotes the strict transform of $ Z $. This proves the assertion
because the image of
the map $H^{2r}(\widetilde Y;\C)\to H^{2r}(X_b;\C)$ is the invariant subspace $I\otimes_{\Q}\C$.

We now observe that, in the hypothesis of the remark, the restriction
$ [Z] \vert_{_{X_b}} $ defines an element in $ J^r(X_b) $ via $ \Psi_{_{AJ}} \, $, for any $ \, X_b \, $ as above.
So, there is a well-defined section (set-theoretical, for the moment) associated to our cycle $ \, Z \, $:
$$
\nu_{_{Z}} \ : \quad \Sigma_d^\circ \ \lra \ \JJ \ .
$$

\rmk{\zfnf}
Denote by $ \, Z = \sum n_i Z_i \, $ the decomposition of $ \, Z \, $ in its irreducible components and
by $ \, L_i \, $ the subspace of $ \, \Sigma_d \, $ parametrizing hypersurfaces containing $ \, Z_i \, $.
For a analytic submanifold $ \, B \subseteq \Sigma_d\setminus(\DD\cup \bigcup L_i) \, $ define
$ \, Z_{i,B}:=\{ (x,b) \in Y_B \, : \, x \in Z_i \cap X_b \} \, $.
Since for any $b\in B$ each $Z_i$ meets properly $X_b$ then $Z_B:=\sum n_i Z_{i,B}$ is a relative cycle of codimension $r$ for the family
$Y_B$, i.e. $Z_B\in \ZZ^r(Y_B)$ and each $Z_{i,B}$ is flat over$B$.
Hence from ([V2], Theorem 7.9) we know that $ \nu_{_{Z}} $ is a normal function on $B$ (i.e. it is holomorphic and horizontal on $B$).
Actually, up to replace $ Z $ with others representatives of its class in $ CH^r(Y) $,
one proves that $\nu_{_{Z}}$ is a normal function on {\it{all}} $\Sigma_d^\circ$.

\rmk{\zfnr}
Besides ordinary cohomology  we will also consider the
intersection cohomology  $IH^*(W;\C)$ for a complex
$m\,$-dimensional irreducible projective variety $W$. Here we
recall some properties which we will use in the sequel: for more
details see ([D2], pg. 154-161). First recall that
$IH^*(W;\C)\cong H^*(W;\C)$ when $W$ is smooth, and that  Poincar\'e
Duality, Lefschetz Hyperplane Theorem and Hard Lefschetz Theorem
still hold for intersection cohomology. Moreover when $W$ has only
isolated singularities then we have $H^{q}(W;\C)\cong IH^{q}(W;\C)$ for
any $q>m$. So, by Poincar\'e Duality, one has  a natural
isomorphism $H_{q}(W;\C)\cong IH^{2m-q}(W;\C)$ for $q>m$. Finally
recall that from Decomposition Theorem it follows that
$IH^*(W;\C)$ is naturally embedded in $H^*(\widetilde W;\C)$ as a
direct summand, where $\widetilde W$ is a desingularization of
$W$. So we have a natural surjection $H_*(\widetilde W;\C)\to
IH^*(W;\C)^{\vee}\cong IH^{2m-*}(W;\C)$. This map identifies with the
push-forward $H_{q}(\widetilde W;\C)\to H_{q}(W;\C)$ when $W$ has
only isolated singularities and for $q>m$.

\prg{2}{Generalized Griffiths' Theorems}

\Thm{\gtha}
Let $ \, Y \subseteq \P^N \, $ denote an irreducible, reduced, projective variety of even dimension $ \, n+1 = 2r \ge 4 \, $,
with isolated singularities.
Let $ \, X \, $ denote the intersection of $ \, Y \, $ with a general hypersurface of degree $ \, d \, $.
Assume the vanishing cohomology $ \, V \, $ is not contained in the middle Hodge component $ \, H^{r,r-1}(X) \oplus H^{r-1,r}(X) \, $.
Then we have
$$
J^r(X)_{0,\,\alg} \quad = \quad 0 \ .
$$
\rm

We want to stress that in the case where $ \, Y \, $ is a smooth
hypersurface in $\P^N$ this theorem is the \lq\lq first part\rq\rq
of Griffiths' Theorem as stated in ([Sh], Theorem 2.2). It has a
strong analogy with N\"oether-Lefschetz Theorem. The present
generalization is obtained by revisiting Shioda's proof ([Sh], pg.
721-722) in view of a result on the monodromy action [DGF]. As for
a generalization of the \lq\lq second part\rq\rq of Griffiths'
Theorem see Theorem {\gthb} below (and compare it with ([V2],
Theorem 8.25)).

\rmk{\rmkd}
$i$) \ For $ \ d \gg 0 \, $, the hypothesis on the vanishing cohomology $ \, V \, $ is automatically satisfied;
\newline
$ii$) \ in the hypothesis of the Theorem, also $ \ J^r(X)_\alg \, $ vanishes in case the homology
space $ \, H_{n+2}(Y;\C) \, $ vanishes
(e.g. when $ \ Y \, \subseteq \, \P^N \ $ is a nodal complete intersection as well as when $ \, Y \, $
is a hypersurface, and not a cone, with at most one ordinary singular point
([D1], (4.6) Corollary. (A) p. 164, and (4.17) Theorem p. 214).

The first statement holds because otherwise the orthogonal complement $ \, I \, $ of $ \, V \, $ would contain $ \, H^{0,\,2r-1}(X) \, $ and
the geometric genus of $ \, X \, $ would be bounded by the geometric genus of $ \, \widetilde Y \, $,
which is impossible for $ \, d \, $ sufficiently large.
As for the second statement, it suffices to observe that the map of Jacobians
$ \ J^r(\widetilde Y) \, \ra \, J^r(X) \ $ is induced by the restriction map
$ \ H^n(\widetilde Y;\, \C) \, \ra \, H^n(X;\, \C) \ $ and that such map factors through $ \, H_{n+2}(Y;\C) \, $,
in fact it is the Poincar\'e dual of the composition
$ \ H_{n+2}(\widetilde Y;\, \C) \, \ra \, H_{n+2}(Y;\, \C) \, \ra \, H_n(X;\, \C) \, $,
where the second map is the Gysin morphism in homology.

\proof{of Theorem $\gtha$}
Let $ \, \{X_t\}_{t\in\P^1} \, $ be a general pencil of degree $ \, d \, $ hypersurface sections and fix a reference point $ \, o \in\P^1\setminus\DD $.
Let $ \, U \subset \P^1\setminus\DD \, $ be a small neighborhood of $ \, o $.
Take a non-zero element $ \, \gamma \in V_o \, $ and extend it by continuous
deformation to $ \, \gamma_t\in V_t \, , \ t \in U \, $.
We have the following dichotomy: either $ \, \gamma_t \in M_t := H^n(X_t;\Q) \cap \left[ H^{r,r-1}(X_t)\oplus H^{r-1,r}(X_t) \right] \, $
for any $ \, t\in\P^1\setminus\DD \, $ and any continuous deformation $ \, \gamma_t \, $ of $ \, \gamma \, , \ $ or the set
$ \, A_{\gamma} \, := \, \{ t\in U \, \vert \, \gamma_t \in M_t \} \, $ is countable.
In fact, by Griffiths' Transversality we know that the condition $ \, \gamma_t \in M_t \, $ is an
analytic condition on $ \, t \in U \, $.
In the first case, the submodule of $ \, V_o \, $ generated by $ \, \gamma \, $ under the action of $ \pi_1(\P^1\setminus\DD,o) \, $ is
contained in $ \, M_o \, $.
Since $ \, V_o \, $ is irreducible ([DGF], Corollary 3.7) then $ \, V_o \subset M_o \, $,
and this is in contrast with our assumption on $ \, V \, $.
So a fortiori $ \, A_{\gamma} \, $ is countable.
Then also $ \ A \, := \, \bigcup_{\gamma\in V_o,\gamma\neq 0} A_{\gamma} \ $ is countable, and for $ t\in U\setminus A $ we have $ M_t\cap V_t=\{0\} $.
Letting $ \, X = X_t $, $ t\in U\setminus A $, on the other hand, as the orthogonal sum $ \ H^n(X;\Q) \, = \, I \perp V \ $
respects the Hodge decomposition, we obtain
$$
I \supseteq M \ , \quad \trm{ equivalently } \quad M \ = \ I \cap \lbrack H^{r,r-1}(X) \oplus H^{r-1,r}(X) \rbrack \ .
$$
Therefore, taking into account that the tangent space to $ J^r(X)_\alg $ at the origin is contained in $ H^{r-1,r}(X) $,
we obtain $ \ J^r(X)_{0,\, \alg} \, = \, 0 \ $ as required.
\qed

\Thm{\gthb} Let $ \ Y \, \subseteq \, \P^N \ $ be an irreducible,
reduced, projective variety of even dimension
$ \, n+1 = 2r \ge 4 \, $, with isolated singularities. Let $ \, Z \in \ZZ^r(Y) \, $ be
a cycle with $ \ cl(Z) \ne 0 \in H_{2r}(Y;\,\C) \, $ and $ \ cl\big([Z]\vert_{X'}\big) = 0 \in H_{2r-2}(X';\,\C) \, $
for some
smooth $ \ d $-degree hypersurface section $ \, X' \, $. Then, for
$ \ d \, \gg \, 0 \ $ and a general $ \ X \ $ smooth $ \ d
$-degree hypersurface section of $ \ Y \, $, the following holds:

\quad {\bf a)\ } $ \ \Psi_{_{AJ}} \big([Z]\vert_X\big) \, $ does not vanish in $ \ J^r(X)_0 \, ; $

\quad {\bf b)\ } $ \ [Z]\vert_X \, $ is not algebraically trivial.

\noindent
Moreover, assuming that the singularities of $ \, Y \, $ are \lq\lq mild" $($see below$)$, then

\quad {\bf a$'$)\ } property {\rm \lq\lq a)"} holds under the weaker hypothesis $ \ d \, \ge \, 3 \, $;

\leftskip15mm \parindent-11mm
\quad {\bf b$'$)\ } property {\rm \lq\lq b)"} holds under the hypotheses $ \ d \, \ge \, 3 \ $ and
that the \lq\lq vanishing cohomology subspace" of $ \ H^{2r-1}(X;\,\C) \ $ is not contained in the
middle cohomology $ \, H^{r,r-1}(X)\oplus H^{r-1,r}(X) \, $.

\leftskip0mm \parindent4mm
\rm

We say that $ \, Y \, $ has {\it mild} singularities if for any $ \, p \in Y \, $ the exceptional divisor of the blow-up of
$ \, Y \, $ at $ \, p \, $ has at worst isolated singularities.

\proof{of Theorem $\gthb$}

\noindent
{\it Step \gthbp1. \ }
We first explain how to deduce the Theorem for $ \ d \, \gg \, 0 \ $ from the following property ($\ssize\bigstar$).

\vskip1mm \leftskip14mm \rightskip0mm \parindent-14mm ($\ssize\bigstar$) \hskip7.5mm
For $ \ d \, \gg \, 0 \, $, \ there exists a $ \ d $-degree hypersurface section $ \ X_o \ $ such that
a Zariski general line $ \ \ell \, \in \, \big\{ \trm{ lines in } \ \Sigma_d \ \trm{ through } X_o \big\} \ $
does not contain a ball $ \ U \, \subseteq \, \ell \setminus \DD $ with
$$
\Psi_{_{AJ}} \, \left([Z]\vert_{X_b}\right) \quad \in \quad T_b \ = \ \image \big( J^r(\widetilde Y) \big) \ \ , \qquad \forall \ b \, \in \, U
$$
(note that this is the vanishing in $ \, J^r(X_b)_0 \, $ condition).

\leftskip0mm \rightskip0mm \parindent4mm

\proof{of Step \gthbp1}
Consider the set $ \ G2 \subseteq \Sigma_d^\circ \ $ of $ \, d $-degree smooth hypersurface sections $ \, X_b \, $ satisfying a).
Since the normal function is analytic and $ \, T \, $ is closed in $ \, \JJ \, $,
the complement of $ \, G2 \, $ in $ \, \Sigma_d^\circ \, $ is analytic.
As a consequence, we have the dichotomy: either $ \, G2 \, $ is empty, or it is the
complement of a proper analytic subset of $ \, \Sigma_d^\circ \, $.
As the first case can be excluded for otherwise ($\ssize\bigstar$) would be contradicted, then $ \, G2 \, $ is dense in $ \, \Sigma_d^\circ \, $.
By Theorem \gtha, which can be applied in view of remark \rmkd,
the set $ \, G1 \, $ of sections satisfying $ \ J^r(X_b)_{0,\,\alg} \, = \, 0 $
contains the complement of a countable union of proper analytic subvarieties of $ \, \Sigma_d \, $.
So, the same holds for the set $ \, G1\cap G2 \, $ as well.
Finally, for $ \, X_b \, $ with $ \ b \ \in \ G1\cap G2 \ $ property b) holds.
In fact, if $ \, [Z]\vert_{X_b} \, $ were algebraically trivial, then
$ \, \Psi_{_{AJ}} \left([Z]\vert_{X_b}\right) \, $ would belong to $ \ J^r(X_b)_{0,\,\alg} \, = \, 0 $
\ (this is because $ \, b \, \in \, G1 \, $), and this would contradict the fact that $ \, b \, \in \, G2 \, $.
\qed

Now the strategy is the following: we introduce a particular hypersurface section $ \, X_o \, $ (step \gthbp2),
we state properties of a Zariski general pencil through $ \, X_o \, $ (step \gthbp3), then we prove property ($\ssize\bigstar$).

\vskip2mm
\noindent
{\it Step \gthbp2. \ }
For $ \ d \, \gg \, 0 \, $, \ there exists a $ \ d $-degree hypersurface section $ \, X_o \, $ of $ \, Y \, $
intersecting properly each components of $ \, Z \, $, \ containing $ \, Y_\sing \, $,
\ and such that $ \, \widetilde X_o \, $ is irreducible, smooth, and very ample,
where $ \, \widetilde X_o \, $ denotes the strict transform of $ \, X_o \, $ in the desingularization
$ \, \widetilde Y \, $ of $ \, Y \, $. Note that, in particular, $ \ X_o \cap Y_\smooth \, $ is smooth.

\proof{of Step \gthbp2}
One can construct $ f: \widetilde Y \to Y $ via a sequence of
blowings-up along smooth centers supported in $ Y_{\sing} $ ([L], Theorem 4.1.3 pg. 241).
Denote by $ g: \widetilde \P \to \P^N $
the projective variety obtained with the same sequence of
blowings-up on $ \P^N $ and observe that $ \widetilde Y\subset \widetilde\P$.
Since the centers are smooth then $ \widetilde \P $ is nonsingular, and its Picard group is generated
by the pull-back of the hyperplane  $ g^*(H) $ and the components
$ E_i\subset \widetilde \P $ of the exceptional divisor
(see [H], Proposition 7.16 and Theorem 8.24, and [F1], Proposition 6.7, (e)).
Therefore for suitable integers $ d \gg 0 $ and $ m_i $ the line bundle $ \OO_{\widetilde \P}\big(dg^*(H)-\sum m_iE_i\big) $ must be very ample.
Let $ \widetilde S_{o} \in \vert \OO_{\widetilde \P}(dg^*(H)-\sum m_iE_i) \vert $ be a general divisor.
Then $ X_o := Y\cap g_*(\widetilde S_{o}) $ verifies all requests.
In fact $ \widetilde X_o $ is equal to $ \widetilde Y \cap \widetilde S_{_0} $, which is irreducible, smooth and very ample.
Moreover $ X_o$ meets properly each component $ Z_i$ of $ Z$ for $ \widetilde X_o $ meets properly $ Z_i\setminus g^{-1}(Y_{\sing}) $.
And $ X_o$ contains $ Y_{\sing}$ because $ \widetilde S_{o} $ is very ample.
\qed

\vskip2mm
\noindent
{\it Step \gthbp3. \ }
A Zariski general pencil $ \, \ell \cong \P^1\, $ through $ \, X_o \, $ satisfies the following properties:

\noindent
$ \bullet \quad $for any $ \ b \, \in \, \ell \ $ the fiber $ \, X_b \, $ of $ \, \ell \, $ meets all components of our cycle $ \, Z \, $
properly, so $ Z_{\ell\backslash\DD} $ is a relative cycle for the family $ Y_{\ell\backslash\DD} $ (cfr remark \zfnf);

\noindent
$ \bullet \quad \ell \backslash \{X_o\} \, $ has finitely many singular fibers, not meeting $ \, Y_\sing \, $ and only having one ordinary double point.

\proof{of Step \gthbp3}
Let $ Z_i $ be an irreducible component of $ Z $.
Since $ \dim(Z_i) = r > 0 $ then $ Z_i $ imposes at least two conditions to $ \Sigma_d $.
Therefore the subspace $ L_i\subset\Sigma_d $ parametrizing hypersurfaces containing $ Z_i $ has codimension at least $ 2 $.
Since $ X_o \not\in L_i $ then a general line in $ \Sigma_d $ passing through $ X_o $ does not meet $ L_i $.
This proves the first property.
\newline
As for the second property, embed $ Y $ via the $ \, d $-Veronese and interpret
$ \Sigma_d $ as the linear series cut by hyperplane sections.
We have $ \, \DD = Y^*\cup p_1^*\cup\dots\cup p^*_h $, where $ Y^* $ denotes the dual
variety of $ Y $, and the $ p^*_i $ denote the dual hyperplanes of the singular points $ p_i $ of $ Y$.
Let $ \ell $ be a general pencil through $ {X_o} $ and let $ \, X_o, \, X_{b_2}, \, \cdots, \, X_{b_k} \, $ be its singular fibers,
namely the fibers corresponding to the points of the intersection $ \, \ell \cap \DD $.
Since $ {X_o}\in \bigcap p^*_i $ then $ \ell $ meets $ \bigcup p^*_i $ only at $ {X_o} $.
Moreover since $ \dim (Y^*)_\sing\leq \dim\Sigma_d-2 $ then $ \ell\cap (Y^*)_\sing\subseteq\{X_o\} $.
Therefore $ \, \{b_2,\,\cdots,\,b_k\} \subset Y^* \setminus \left[(Y^*)_\sing\cup\bigcup p_i^*\right]$.
We may assume that $ \{b_2,\,\cdots,\,b_k\}$ is not empty otherwise the proof of the second property ends here.
In particular we may assume that $ Y^*$ has no dual defect, and that it is not a cone with vertex $ X_o$.
Now we claim that, for any $ \, j \, = \, 2,\,\cdots,\,k$,
the pencil $ \ell $ meets $ Y^*$ transversally at $ b_j $ (this means that $ X_{b_j} $ only has one ordinary double point).
Suppose not. Then the general line through $ X_o $ is tangent to $ Y^*\setminus (Y^*)_\sing $ at some point.
It follows that if we consider the projection (possibly internal) of $ Y^* $ from
$ X_o $, then the image of $ R \setminus \{X_o\} $ is equal to the image of $ Y^*\setminus\{X_o\} $,
where $ R := \{b\in Y^*\,:\,X_o\in T_{b,Y^*}\} $ denotes the ramification locus of the projection.
Since the image of $ Y^*\setminus\{X_o\} $ has the same dimension of $ Y^* $
because $ Y^*$ is not a cone with vertex $ X_o $, then the ramification locus is all $ Y^* $.
Therefore for a general $ b\in Y^* $ we have $ X_o \in T_{b,Y^*} $.
Now let $ J \subset \Sigma_d $ be the cone with vertex $ X_o $ and basis $ Y^* $ (i.e. the embedded join of $ X_o $ and $ Y^* $).
By Terracini's Lemma ([FOV], Proposition 4.3.2.) we know that for a general $ b\in Y^* $ and a general $ c \in \overline{bX_o} $
the tangent space to $ J $ at $ c $ is spanned by $ X_o $ and the tangent space to $ Y^* $ at $ b $.
Since $ X_o\in T_{b,Y^*} $ then $ \dim J = \dim Y^* $, i.e. $ J=Y^* $.
Hence $ Y^* $ is a cone with vertex $ X_o $, and this is in contrast with previous assumption.
\qed

Now, for a pencil $ \, \ell \cong \P^1 \, $ as in step \gthbp3, let $ \ B \, = \, \P^1 \setminus \DD \ $
denote the set of smooth sections
and consider the natural completion  of the family $ \ \pi : \, Y_B \ra B \ $ (cfr \S 1)
$$
\CD
Y_B                        @. \quad \hookrightarrow \quad @.  Y_\ell              \\
\pi \ \downarrow \hskip4mm @.                             @. \hskip4mm \downarrow \ \pi \\
B                          @. \hookrightarrow             @. \P^1
\endCD
$$
where $ \,  Y_\ell \, $ is the blow-up of $ \,  Y \, $ along the base locus of the pencil.
Finally, consider $ \, Z_B \, \in \, \ZZ^{r}(Y_B) \, $ as introduced in \S1, remark  \zfnf, and its class
$ \ cl(Z_B)^\vee \, \in \, H^{2r}(Y_B; \, \Z) \, $.

We are now ready to prove property ($\ssize\bigstar$). For this, we proceed by contradiction.

\noindent
{\it Step \gthbp4. \ }
Let $ \, \ell \, $ be a general pencil as in step \gthbp3.
If there exists a ball $ \, U \subseteq B \, $ such that $ \, \nu_{_Z}(b) \in T_b \, $ for all $ \, b \in U $, then there exists an element
$ \, \xi \in J^r(\widetilde Y) \, $ such that
$ \ i_b^*\,(\xi) \, = \, \nu_{_Z}(b) \, , \ \forall \, b \in B $.

\proof{of Step \gthbp4}
First, we introduce a natural commutative diagram
$$
\CD
@. @. J^r(\widetilde Y) \\
@. @. @VV \, {i_b^*} V \\
{ [Z]\vert_{X_b} \quad \in } @. CH^r(X_b)_{\hom} @> \Psi_{_{AJ}} >> {J^r(X_b)} @. \quad \ni \quad \nu_{_Z}(b) \\
@. @VV V @VV \, {i_b}_\star V \\
{ i_{b \,*}}({ [Z]\vert_{X_b}}) { \quad \in } \quad @. CH^{r+1}(\widetilde Y)_{\hom} @> \Psi_{_{AJ}} >> J^{r+1}(\widetilde Y)
\endCD
$$
where the left vertical map denotes the push-forward $ \ i_{b \,*} \ $ and the map $ \, {i_b}_\star \, $ is the map induced
by the map $ \ {i_b}_\star \, : \ H^{2r-1}(X_b;\,\Z) \, \ra \, H^{2r+1}(\widetilde Y;\,\Z) \ $
(which is the Gysin morphism in cohomology associated to $ \ i_b \, : \ X_b \, \ra \, \widetilde Y \, $).
We now observe the following properties:
\newline
$ i) \ $ the vertical composition $ \ {i_b}_\star \, \circ \, i_b^* \ $ can be interpreted as the map induced
by the map
$ \ H^{2r-1}(\widetilde Y;\,\Z) \, \ra \, H^{2r+1}(\widetilde Y;\,\Z) $ \ given by the
cap product in homology (modulo Poincar\'e duality) with the class of the \lq\lq divisor
$ \, X_b \, $" (in particular it does not depend on $ \, b \, $);
\newline
$ ii) \ $ as the diagram above commutes and $ \ i_{b \,*}({ [Z]\vert_{X_b}}) \, \in \, CH^{r+1}(\widetilde Y)_{\hom} \ $ does not depend on $ \, b \, $,
the image $ \ {i_b}_\star \circ \nu_{_Z}(b) \ $ does not depend on $ \, b \, $ as well.
\newline
We now conclude the proof of step \gthbp4 under the assumption that the restriction $ \ {i_b}_\star \vert_{\image (i_b^*)} \ $ is an isogeny on its image.
Let us now go back to our exact sequence
$ \ 0\ra K \times B \, \ra \, J^r(\widetilde Y) \times B \, \ra \, T \, \ra \, 0 \ $
of families over $ \, B \, $.
Working modulo the identification $ \ {J^r(\widetilde Y) \over K} \times B \ \cong \ T \ $ of fibrations over $ \, B \, $, we replace
$$
J^r(\widetilde Y) \times B \quad \blr{\lra}{i^*} \quad T \quad \blr{\lra}{i_\star} \quad J^{r+1}(\widetilde Y) \times B
$$
with
$$
J^r(\widetilde Y) \times B \quad \lra \quad {J^r(\widetilde Y) \over K} \times B \quad \blr{\lra}{\hat \iota_\star} \quad J^{r+1}(\widetilde Y) \times B
$$
where $ i_\star $ denotes the map of Jacobian fibrations that globalizes $ {i_b}_\star $.
Now, our normal function $ \, \nu_{_Z} \, $ takes values to $ \, T \, $  and therefore induces a analytic section
$ \ \hat\nu \, : \ B \ \ra \ {J^r(\widetilde Y) \over K} \times B $.
\ As the image $ \ {i_b}_\star \circ \nu_{_Z} (b) \ $ does not depend on $ \, b $, \ the composition
$ \, \hat \iota_\star \circ \hat\nu \, $ is constant
(as a matter of language, a section of a trivial fibration is said to be {\it {constant}} if its image is the graph of a constant function).
By our assumption stating that the restriction $ \ {i_b}_\star \vert_{\image (i_b^*)} \ $ is an isogeny on its image,
the kernel of $ \, \hat \iota_\star \ $ is discrete
\ (by definition, the kernel of a map of fibrations is the inverse image of the zero section, and it is said to be
{\it{discrete}} if
its restriction to any fiber is a discrete subset).
As a consequence, the inverse image of $ \ \hat\nu \ $ itself must be constant.
\newline
We are left to prove that the restriction $ \ {i_b}_\star \vert_{\image (i_b^*)} \ $ is an isogeny on its image.
First we may note that in the smooth case $ \, X_b \, $ is ample on $ \, Y \, = \, \widetilde Y \, $ and
$ \ {i_b}_\star \, \circ \, i_b^* \ $ is an isogeny by the Hard-Lefschetz theorem (and the interpretation we have given above).
\newline
Going back to our situation, we set $ \ i_b , \ j_b, \ f \ $ as in the diagram
$$
\CD
@. \widetilde Y \\
{\ssize i_b} \nearrow \quad @. \downarrow @. {\ssize f} \\
X_b \quad \blr{\lra}{j_b} \quad @. Y
\endCD
$$
and we investigate more closely the map $ \ {i_b}_\star \, \circ \, i_b^* \, $.
As we already said, the map $ \ {i_b}_\star \, \circ \, i_b^* \ $ descends from the map
$ H^{2r-1}(\widetilde Y;\,\Z) \ \lra \ H^{2r+1}(\widetilde Y;\,\Z) \ $ which, passing to complex coefficients,
turns out to be equal to the composition of all maps at the first row in the diagram below (cfr remark \zfnr \ on intersection cohomology):
$$
\CD
H^{2r-1}(\widetilde Y;\,\C) @. \ \cong \ @. H_{2r+1}(\widetilde Y;\,\C) @. \ \blr{\lra}{f_*} \ @. H_{2r+1}(Y;\,\C) @.
\ \cong \ @. IH^{2r-1}(Y;\,\C) @. \ \cong^{^{(1)}} \ @. IH^{2r+1}(Y;\,\C) @. \ \subseteq \ @. H^{2r+1}(\widetilde Y;\,\C) \\
\downarrow {\ssize i_b^* } @. @. @. @. \downarrow {\ssize j_b^\star } @. @. @. @. @. @. \updownarrow {\ssize \cong} \\
H^{2r-1}(X_b;\,\C) @. @. \cong @. @. H_{2r-1}(X_b;\,\C) @. @. @. \blr{\lra}{i_{b \,*}} @. @. @. H_{2r-1}(\widetilde Y;\C).
\endCD
$$
Here $ \, i_b^* \, $ is the natural pull-back and $ \, j_b^\star \, $ is the Gysin morphism in homology
and the isomorphism (1) is the Hard-Lefschetz isomorphism.
Note that the map of Jacobians $ \ i_b^* \, : \ J^r(\widetilde Y) \, \ra \, J^r(X_b) \ $ descends from the first one
and $ \ {i_b}_\star \ $ descends from the composition of the bottom row with the right \lq\lq vertical" isomorphism (Poincar\'e duality).
Finally, an element in $ \, H^{2r-1}(\widetilde Y;\,\C) \, $ that maps to zero in $ \, H^{2r+1}(\widetilde Y;\,\C) \, $
must vanish in $ \, H_{2r+1}(Y;\,\C) \, $ and therefore must vanish in $ \, H^{2r-1}(X_b;\,\C) \, $.
\newline
This proves that the restriction of the bottom row to the image of $ \, i_b^* \, $ is injective, and therefore that the restriction
$ \ {i_b}_\star \vert_{\image (i_b^*)} \ $ is an isogeny on its image as required.
\qed

\newpage

\noindent {\it Step \gthbp5. \ }
The cohomology class of $ Z_B $ in $ H^{2r}\big(Y_B;\,\C\big) $ is zero.

\proof{of Step \gthbp5}
First recall the exact sequence defining the Jacobian fibration:
$ 0\to R^{2r-1}\pi_*\Z\slash{\text{torsion}}\to \HH^{2r-1}\slash F^r\HH^{2r-1}\to \JJ \to 0, $
where $\HH^{2r-1}:=R^{2r-1}\pi_*\Z\otimes \OO_B$ ([V2], (8.12), p. 230).
From previous step we know there exists $\xi\in J^r(\widetilde Y)$ such that $ i_b^*(\xi) = \nu_{_Z}(b) $ for any $ b \in B $.
Let $ \xi'$ be a lift of $ \xi$ in $ H^{2r-1}(\widetilde Y;\,\C)$.
Then $ \{i_b^*(\xi')\otimes 1\}_{b\in B}\in H^0(B, \HH^{2r-1}\slash F^r\HH^{2r-1}) $
is a global section whose image in $ H^0(B, \JJ)$ is equal to $ \nu_{_Z} $
(here, $ i_b^* $ is the pull-back $ H^{2r-1}(\widetilde Y,\C)\to H^{2r-1}(X_b,\C) $).
Therefore the image of $ \nu_{_Z}$ in $ H^1(B,R^{2r-1}\pi_*\C) $ vanishes.
On the other hand the Leray filtration of $ Y_B\to B$ induces a natural map
$ \ker(H^{2r}(Y_B;\C)\to H^{2r}(X_b;\C))\to H^1(B,R^{2r-1}\pi_*\C) $,
which is an isomorphism because $ B\subset \P^1 $ (compare with [V2], Theorem 8.21 and proof of (ii)).
And under this identification one knows that the image of $ \nu_{_Z} $ corresponds to $ Z_B$ ([V2], Lemma 8.20).
\qed

Remembering our strategy, to find a contradiction, hence to conclude the proof of the Theorem for $ \, d \gg 0 \, $,
it suffices to prove step \gthbp6 below
(here we make a strong use of the fact that $ \ell $ is as in step \gthbp3, a assumption that was not necessary for steps \gthbp4 and \gthbp5).

\noindent
{\it Step \gthbp6. \ }
Keep our previous notation, in particular $ \, \ell \, $ as in step \gthbp3 and $ \, B \, = \, \P^1 \setminus \DD \, $. Then,
$$
cl(Z_B)^\vee \ = \ 0 \ \in \ H^{2r}\big(Y_B;\,\C\big) \qquad \Longrightarrow \qquad cl(Z) \ = \ 0 \ \in \ H_{2r}(Y;\,\C) \, .
$$

\proof{of Step \gthbp6}
Let $Y_{\ell}$ be the blow-up of $Y$
along the base locus of $\ell$,  let  $\pi:Y_{\ell}\to \ell\cong\P^1$ be the natural map, and
consider the exact sequence of the pair \ ($ Y_\ell , \, \pi^{-1}(\DD)$)\,, \ maps and elements as indicated in the diagram below:
$$
\CD
@. H_{2r}(Y;\,\C) @. H^{2r}(Y_B;\,\C) @. \ \ni \ @. cl(Z_B)^\vee @. \quad { ( \ = \ 0)} \\
\quad {}^{^{\oplus {j_{b_i *}}}} \nearrow @. @VV \lambda^\star V @VV \cong V @. \updownarrow \\
\oplus H_{2r}(X_{b_i};\,\C) \quad @> \varpi >> \quad H_{2r}(Y_\ell;\,\C) \quad @>\rho >> \quad H_{2r}(Y_\ell,\,\pi^{-1}(\DD);\,\C) @. \ @. \\
\hskip2mm { \sum \xi_i} \hskip13mm \ms \hskip-15mm @. \hskip3mm { \lambda^\star(cl(Z)) } \hskip11mm \ms \hskip-14mm @. @. @. { \rho(\lambda^\star(cl(Z)))} @. \quad { ( \ = \ 0)}
\endCD
$$
where,
$ i) $ \ $ \{b_0 = o, \, b_1,\, ...,\,b_k \} \ = \ \DD \cap \ell \ = \ \ell \backslash B \ $ \ denotes the discriminant locus of the pencil;
$ ii) $ \ $ j_{b_i} \, : \ X_{b_i} \, \hookrightarrow \, Y \ $ \ denotes the natural inclusion;
$ iii) $ \ $ \lambda \, : \ Y_\ell \, \ra \, Y \ $ \ denotes the natural projection, and so
$ \ \lambda^\star \, : \ H_{2r}(Y;\,\C) \, \ra \, H_{2r}(Y_\ell;\,\C) \ $ is the Gysin morphism in homology (notice that
$\lambda$ is a l.c.i. morphism  because the base locus of $\ell$ is contained in $Y_{\smooth}$);
$ iv) $ \ we work under the natural identification $ \ H_{2r}(\pi^{-1}(\DD);\,\C) \ \cong \ \oplus H_{2r}(X_{b_i};\,\C) \ $
induced by the disjoint union decomposition \ $ \pi^{-1}(\DD) \, = \, \cup^\cdot X_{b_i} \, $,
\ note that $ \ \varpi \, = \, \oplus \iota_{b_i *} \ $ where $ \ \iota_{b_i} \, : \ X_{b_i} \, \hookrightarrow \, Y_\ell \ $ is the natural inclusion.

\vskip3mm \leftskip0mm \parindent0mm
We claim that, as $ \ cl(Z_B)^\vee \, = \, 0 \ $ by hypothesis,

$ (\gthbp6.1) $ \qquad there exists an element $ \ \sum \xi_i \ $ as indicated in the diagram and satisfying
$ \quad \sum \, j_{b_i *} ( \xi_i) \ = \ cl(Z) \ $.

The proof of (\gthbp6.1)
involves two statements, the first one is that $ \ cl(Z_B)^\vee \ $ and $ \ \rho(\lambda^\star(cl(Z))) \ $ correspond to each other under Lefschetz Duality
[Sp] (the vertical isomorphism at the right of the diagram), and this is clear.
So, as a consequence of the exactness of the pair sequence there exists an element $ \ \sum \xi_i \ $
mapping to $ \ \lambda^\star(cl(Z)) \ $ via
$ \ \varpi \ $ above.
The second one is the chain of equalities
$$
\left[ \, \oplus \, {j_{b_i *}}\right] \left(\sum \xi_i \right) \quad = \quad
\lambda_* \, \varpi \, \left(\sum \xi_i \right) \quad = \quad \lambda_* ( \lambda^\star(cl(Z))) \quad = \quad cl(Z)
$$
where, the first equality is trivial, the second one follows by the definition of $ \ \sum \xi_i \, $,
the third one is the non-trivial one and follows by ([F1], Proposition 6.7, $ (b)$).
This concludes the proof of claim (\gthbp6.1).

\vskip1mm
Now, for any $ \ b \, \in \, \P^1 \, $, \ we consider
the Gysin map $ \ j_b^\star \ $ associated to the inclusion $ j_b \, : \ X_b \, \hookrightarrow \, Y \ $ and the diagram
$$
\CD
H_{2r+2}(Y;\,\C) @. \quad @> \gamma >> \quad H_{2r}(Y;\,\C) \\
@V j_b^\star VV @. \hskip-10mm \nearrow {\ssize {j_{b \,*}}} \hskip10mm \\
H_{2r}(X_b;\,\C) @.
\endCD
\tag{\gthbp6.2}
$$
where $ \, \gamma \, $ denotes the cap-product with
$cl(X)^\vee\in H^2(Y; \, \C)$
(namely $ \, \gamma = \cap \,cl(X)^\vee \, $, where $ \, X \, $ is any $ \, d $-degree hypersurface section of $ \, Y \, $).
This diagram is commutative and therefore, in particular, the
composition $ \, j_{b \,*} \circ j_b^\star \, $ does not depend on $ \, b $.
We now show the following
\vskip3mm
$ (\gthbp6.3) $ \hfill for any $ \ i \, \in \, \{0,\,...,\,k \} \ $ the Gysin map $ \ j_{b_i}^\star \ $ is surjective. \hfill \
\vskip3mm

To prove (\gthbp6.3),
we first examine the case $ i=0$. Consider the following natural commutative diagram:
$$
\CD
@. H_{2r+2}(\widetilde Y;\C) @> >> H_{2r}(\widetilde X_o;\C) @. \\
@. @VV V @VV V \\
@. H_{2r+2}(Y;\C) @> j_o^\star
>> H_{2r}(X_o;\C),
\endCD
$$
where the horizontal maps denote Gysin maps, and the vertical ones
denote push-forward. Since $ \widetilde X_o$ is a smooth
and very ample divisor on $ \widetilde Y$ (see step \gthbp2) then the upper
horizontal map is an isomorphism in view of Poincar\'e Duality and
Hyperplane Lefschetz Theorem. Moreover the vertical map at the
right is surjective: this follows by the Decomposition Theorem because
$ f_{|\widetilde X_o}:\widetilde X_o\to X_o$ is a desingularization of
$ X_o$ (see remark {\zfnr}).
Then the commutativity of the diagram implies that $ j_o^\star $ is surjective.

Next assume $ i\in\{1,\dots,k\}$, i.e. assume that $
X:=X_{b_i}$ only has one ordinary double point at a smooth point
$ p$ of $ Y$. Denote by $ B_p(Y)$ and $
B_p(\P^N)$ the blowing-up of $ Y$ and $ \P^N$ at $
p$, by $ E_Y$ and $ E_{\P}$ the exceptional divisors,
and by $ \phi:B_p(Y)\to Y$ and $ \psi:B_p(\P^N)\to\P^N$
the natural projections. Let $ \widetilde X$ be the strict
transform of $ X$ in $ B_p(Y)$. Then $
\widetilde X$ is a smooth Cartier divisor on $ B_p(Y)$, and
also very ample on $ B_p(Y)$: in fact $ \OO_{B_p(Y)}(\widetilde X)\cong \OO_{B_p(Y)}(\phi^*(dH_Y)-2E_Y)$
is the restriction on $ B_p(Y)$ of
$ \OO_{B_p(\P^N)}(\psi^*(dH_{\P})-2E_{\P})$, which is very
ample on $ B_p(\P^N)$ (here we denote by $ H_{\P}$ the
hyperplane in $\P^N$, and by $ H_Y$ its restriction to
$ Y\subset\P^N$). Now, as before, consider the natural
commutative diagram:
$$
\CD
@. H_{2r+2}(B_p(Y);\C) @> >> H_{2r}(\widetilde X;\C) @. \\
@. @VV V @VV V \\
@. H_{2r+2}(Y;\C) @>>> H_{2r}(X;\C).
\endCD
$$
As before the right vertical map is surjective because $
\widetilde X$ is a desingularization of $ X$. Moreover,
since $ \widetilde X$ is smooth and very ample on $
B_p(Y)$ then by Lefschetz Hyperplane Theorem for "intersection
cohomology" the restriction map $ IH^{2r-2}(B_p(Y);\C)\to
IH^{2r-2}(\widetilde X;\C)$ is an isomorphism (remark {\zfnr}). But
$ IH^{2r-2}(B_p(Y);\C)\cong H_{2r+2}(B_p(Y);\C)$ because $
B_p(Y)$ only has isolated singularities, and $
IH^{2r-2}(\widetilde X;\C)\cong H^{2r-2}(\widetilde X;\C)\cong
H_{2r}(\widetilde X;\C)$ because $ \widetilde X$ is smooth
and by Poincar\'e Duality. Putting together it follows that
$ H_{2r+2}(B_p(Y);\C)\to H_{2r}(\widetilde X;\C)$ is
surjective. The commutativity of the diagram shows then that also
$ H_{2r+2}(Y;\C)\to H_{2r}(X;\C)$ is surjective. This
concludes the proof of $ (\gthbp6.3)$.

\vskip3mm
We now conclude the proof of step \gthbp6.
By (\gthbp6.3), for any $ \ i \, \in \, \{ 0,\,...,\,k \} \ $ there exists an element
$ \, \eta_i \, $ satisfying $ \ j_{b_i}^\star(\eta_i) \, = \, \xi_i \ $. \ Thus, we have
$$
cl(Z) \quad = \quad \sum \, j_{b_i *}( \xi_i) \quad = \quad \sum j_{b_i *}\left(j_{b_i}^\star(\eta_i)\right) \quad = \quad
\sum \gamma(\eta_i) \quad = \quad \gamma \left(\sum \eta_i \right)
$$
where the first equality is (\gthbp6.1), the second one is obtained substituting $ \, \xi_i \, $ with its expression
$ \, j_{b_i}^*(\eta_i) \ $, the third one follows by the commutativity of diagram (\gthbp6.2)
and the last one is trivial.

Going back to our diagram (\gthbp6.2), the equality $ \ cl(Z) \, = \, \gamma \left(\sum \eta_i \right) \ $ shows that $ \, cl(Z) \, $
belongs to the image of the push-forward $ \ j_{b \,*} \, : \ H_{2r}(X_b;\,\C) \, \ra \, H_{2r}(Y;\,\C) \, $,
i.e. there exists $ \ \alpha \, \in \, H_{2r}(X_b;\,\C) \ $ satisfying $ \ cl(Z) \, = \, j_{b \,*}(\alpha) \, $.
On the other hand, choosing a smooth $ X_b $, we might consider the composition
$$
\CD
H_{2r}(X_b;\,\C) \quad @>j_{b \,*}>> \quad H_{2r}(Y;\,\C) \quad @>j_b^\star>> \quad H_{2r-2}(X_b;\,\C) \\
\alpha \hskip14mm \ms \hskip-16mm @. cl(Z) \hskip15mm \ms \hskip-17mm @. cl([Z]\vert_{X_b}) \, = \, 0
\endCD
$$
and observing it is injective by the Hard-Lefschetz theorem, we deduce that $ \ cl(Z) \, = \, 0 \ $ as required.
\qed

\vskip8mm
So far we have proved {\bf a)} and {\bf b)}. We now prove {\bf a$'$)} and {\bf b$'$)}.

The previous proof can be adapted with minor changes. First, as before, we have that

\def\bis{${}_{\operatorname{bis}}$}

\vskip4mm
\noindent
{\it Step \gthbp1\bis. \ }
The theorem follows from the following property (${\ssize\bigstar}'$).

\vskip1mm \leftskip14mm \rightskip0mm \parindent-14mm (${\ssize\bigstar}'$) \hskip7.5mm
For $ \ d \, \ge \, 3 \, $, the Zariski general line $ \ \ell \, \in \, \big\{ \trm{ lines in } \ \Sigma_d \ \big\} \ $
does not contain a ball $ \ U \, \subseteq \, B $ \ with
$ \ \Psi_{_{AJ}} \, \left([Z]\vert_{X_b}\right) \quad \in \quad \image \big( J^r(\widetilde Y) \big) \ \ , \qquad \forall \ b \, \in \, U \, $.

\leftskip0mm \rightskip0mm \parindent4mm

\noindent
{\it Step \gthbp3\bis. \ }
A Zariski general pencil $ \ \ell \ $ satisfies the following properties:

\noindent
$ \bullet \quad $for any $ \ b \, \in \, \P^1 \ $ the fiber $ \, X_b \, $ of $ \, \ell \, $ meets all components of our cycle $ \, Z \, $ properly;

\leftskip5mm \parindent-5mm
$ \bullet \quad \ell \ $ has finitely many singular fibers: those not meeting $ \, Y_\sing \, $ only having one ordinary double point,
and those meeting $ \, Y_\sing \, $ only at one point and generically
(so, any of such sections is singular only at one point and its singularity is a
general section of a \lq\lq mild" singularity).
\vskip0mm
\leftskip0mm \parindent0mm

\vskip1mm
\noindent
{\it Step \gthbp4\bis. \ }
If there exists a ball $ \ U \, \subseteq \, B \ $ such that $ \ \nu_{_Z}(b) \, \in \, T_b \ $ for all
$ \ b \, \in \, U \, $, \ then there exists an element
$ \ \xi \, \in \, J^r(\widetilde Y) \ $ such that
$i_b^*\,(\xi) =  \nu_{_Z}(b)$ for any $b \ \in \ B$.

\vskip1mm
\noindent
{\it Step \gthbp5\bis. \ }
The cohomology class of $ \, Z_B \, $ in $ \ H^{2r}\big(Y_B;\,\C\big) \ $ is zero.

\vskip1mm
\noindent
{\it Step \gthbp6\bis. \ }
We consider a pencil $ \, \ell \, $ as in step \gthbp3\bis. Then,
$$
cl(Z_B)^\vee \ = \ 0 \ \in \ H^{2r}\big(Y_B;\,\C\big) \qquad \Longrightarrow \qquad cl(Z) \ = \ 0 \ \in \ H_{2r}(Y;\,\C). \,
$$

The proof of steps \gthbp1\bis, \gthbp3\bis, \gthbp4\bis, \gthbp5\bis
are those of the analogous steps from the proof of {\bf a)} and {\bf b)}. We now prove step \gthbp6\bis.

\proof{of Step \gthbp6\bis}
Using the same argument as in the proof of step \gthbp6
before, in order to prove step \gthbp6\bis \ one reduces to prove that for
any $ \ i \, \in \, \{ 1,\,...,\,k \} \ $ the Gysin map $ \, j_{b_i}^\star \, $ is surjective, where $ \{b_1,\dots,b_k\} $
denotes the discriminant locus of the pencil
(compare with \gthbp6.3 in the proof of step \gthbp6). Put $ X:=X_{b_i} $.

When $ X$ only has one ordinary double point at a smooth
point $ p$ of $ Y$ then the same argument we used in
this case in the proof of (\gthbp6.3) applies (here the assumption
$ d\geq 3$ allows us to say that $ \OO_{B_p(\P^N)}(\psi^*(dH_{\P})-2E_{\P})$ is very ample).

It remains to examine the case where $ X $ is a general
hypersurface section through a singular point $ p$ of
$ Y$. Since the singularity is "mild" the analysis is quite
similar to the previous case. In fact, denote by $ B_p(Y)$
and $ B_p(\P^N)$ the blowing-up of $ Y$ and $ \P^N$
at $ p$, by $ E_Y$ and $ E_{\P}$ the exceptional
divisors, and by $ \phi:B_p(Y)\to Y$ and $
\psi:B_p(\P^N)\to\P^N$ the natural projections. Let $ \widetilde X$ be the strict transform of $ X$ in $ B_p(Y)$. Since
the singuar point $ p$ is "mild" then $ B_p(Y)$ still
has only isolated singularities and $ \widetilde X$ is a
smooth Cartier divisor on $ B_p(Y)$. Moreover
$\widetilde X$ is also very ample on $ B_p(Y)$ because
$ \OO_{B_p(Y)}(\widetilde X)\cong \OO_{B_p(Y)}(\phi^*(dH_Y)-E_Y)$, and this line bundle is the
restriction on $ B_p(Y)$ of $ \OO_{B_p(\P^N)}(\psi^*(dH_{\P})-E_{\P})$, which is very ample on
$ B_p(\P^N)$ because $ d\geq 2$. At this point one may
prove that $ j_{b_i}^\star $ is surjective exactly as in the
"tangential" case. This concludes the proof of step \gthbp6\bis.
\qed

\enddocument